\documentclass[11pt]{amsart}
\usepackage{graphicx} 
\usepackage{amsmath}
\usepackage{amsthm}
\usepackage{amssymb}
\usepackage{amsfonts}
\usepackage{dsfont}
\usepackage{cite}
\usepackage{hyperref}
\newtheorem{theorem}{Theorem}
\newtheorem{lemma}{Lemma}
\newtheorem{corollary}{Corollary}
\newtheorem{question}{Question}
\newtheorem*{theorem*}{Main Theorem}
\begin{document}

\title{Positive Codegree Density of Suspensions of Graphs}

\author[Haiwen Luo]{Haiwen Luo\textsuperscript{*}}

\thanks{\textsuperscript{*}Qiuzhen College, Tsinghua University, Beijing, China. 
Email: \href{mailto:luohw25@mails.tsinghua.edu.cn}
{\nolinkurl{luohw25@mails.tsinghua.edu.cn}}.}



\begin{abstract}
We determine the exact positive codegree density of the $r$-uniform suspensions $S^r(G)$ of graphs $G$ when $\chi(G)=t$ and $t-2$ is a prime power. The value is $\gamma^+(S^{r}(G)) = \frac{(t-3)(t-2)^{r-1}}{(t-2)^r-1}$. The case of suspensions of $K_4$ gives the positive codegree density of daisies, answering a question of Balogh, Halfpap, Lidick\'y and Palmer. The lower bound construction follows from Ellis, Ivan and Leader, while we establish the matching upper bound based on a weighted graph inequality.
\end{abstract}

\maketitle

\section{Introduction}
Hypergraph Tur\'an problems deal with extremal hypergraphs without a given hypergraph or hypergraph family. For a general introduction, see \cite{htp}. The hypergraph Tur\'an density is defined to be 
\[
\pi(G) = \lim_{n\to\infty}\frac{\mathrm{ex}(n,G)}{\binom{n}{r}},
\]
meaning the maximum possible asymptotic density of a $G$-free hypergraph.

Studying Tur\'an density is an active area of research, but until now we only know the exact Tur\'an densities of relatively few hypergraphs. In contrast, the problem of determining Tur\'an density for graphs was already solved by Erd\H{o}s, Stone and Simonovits in \cite{esi_g,est_g}. This shows that determining Tur\'an densities for hypergraphs is more difficult than for graphs.

Besides Tur\'an density, one can also study extremal parameters defined by degree. This results in the definition of codegree density, which was established by Mubayi and Zhao in \cite{co_d}, and positive codegree density, which is the main focus of our article.

For an $r$-uniform hypergraph $H$, we use $d_H(S)$ to denote the number of edges containing $S$. We define the minimum positive codegree of $H$ to be
\[
\delta^+(H)  =\min_{|F| = r-1\text{ is contained in some edge} }d_H(F).
\]
For a hypergraph $G$, the positive codegree Tur\'an number $\mathrm{coex}^+(n,G)$ is defined as the maximum of $\delta^+(H)$ over all $G$-free $n$-vertex hypergraphs. The positive codegree (Tur\'an) density is defined as
\[
\gamma^+(G)=\lim_{n\to\infty}\frac{\mathrm{coex}^+(n,G)}{n}.
\]
The notion of positive codegree density was introduced in the foundational article \cite{pco_b}, where Halfpap, Lemons and Palmer proved that the limit exists and studied the properties of positive codegree density.

The positive codegree density of hypergraphs has attracted increasing attention over the last few years \cite{AES_gt,np_pco,pco_C5,pcdj}. Nevertheless, exact values have been determined for only a limited range of hypergraphs. Known examples include several small $3$-uniform hypergraphs \cite{pcdj,pco_b}, the $r$-triangle $T^r=S^r(K_3)$ \cite{pcdj}, certain $3$-uniform tight cycles and tight cycles with one edge deleted \cite{pcdj,pco_C5}, and generalized triangles \cite{AES_gt}. A closely related line of research is studying the positive codegree thresholds for the existence of spanning structures, including perfect matchings \cite{pco_M,pco_sp}, Hamilton cycles \cite{pco_H,pco_sp,sco_H} and spanning spheres \cite{sp_sp}.

Our article focuses on the positive codegree density of daisies, and more generally, suspensions of graphs. For a graph $G$ (or $k$-uniform hypergraph $G$, $k<r$), let $S$ be a set of $r-2$ ($r-k$) new vertices
disjoint from $V(G)$. The $r$-uniform suspension $S^r(G)$ is an $r$-uniform hypergraph with vertex set $V(G)\cup S$ and edge set
\[
E(S^r(G))=\{e\cup S:e\in E(G)\}.
\]
In particular, we define the daisy $\mathcal D_r$ to be the $r$-uniform suspension $S^r(K_4)$. More generally, we define the generalized daisy $\mathcal D_{r,t}$ to be the $r$-uniform suspension $S^r(K_t)$.

The study of daisies was popularized by Bollob\'as, Leader and Malvenuto in \cite{cj_dt}, where they proposed several conjectures regarding Tur\'an problems of daisies. One of the conjectures was that daisies have Tur\'an density tending to $0$ as $r$ tends to infinity. The conjecture was disproved by Ellis, Ivan and Leader in \cite{da_co}, where they used a construction based on linear algebra and finite fields. When $t-2$ is a prime power, the positive codegree density of their construction is $\frac{(t-3)(t-2)^{r-1}}{(t-2)^r-1}$, and we prove that it is optimal. They also disproved some conjectures regarding Tur\'an problems in the hypercube, showing the close relationship between the two fields. The study of graph suspensions is also related to coding theory. Alon \cite{e_codes} used results regarding Tur\'an densities of graph suspensions (which he called graph augmentations) to derive a new lower bound for erasure list-decodable codes. 

Balogh, Halfpap, Lidick\'y and Palmer \cite{pcdj} determined the positive codegree density for several hypergraphs and hypergraph families. In particular, they proved that the daisy $\mathcal D_3$ (which they denoted by $J_4$) has positive codegree density $\frac47$, attained asymptotically by balanced blow-ups of the complement of the Fano plane. Their proof combines flag algebra computations and a stability argument. In \cite[Question 6.6]{pcdj} the authors posed the problem of determining the positive codegree density of all daisies $\mathcal D_r$.

Recently, Hou, Hou, Liu, Zeng and Zhang proved that
\[
\gamma^+(\mathcal D_{r,t})\le 1-\frac1{t-1}-\frac1{36(t-1)^2}
\]
for $r\ge 3$ in \cite{F_link}, determining the asymptotic positive codegree density up to $\Theta(t^{-2})$ when $t-2$ is a prime power.

In this article, the author determines the positive codegree density of all $r$-suspensions of graphs $G$ when $\chi(G)=t$ and $t-2$ is a prime power. For $t=4$, this answers the aforementioned problem on positive codegree density of daisies. The following is the main result of the paper.

\begin{theorem*}
Let $r \ge 3$ be arbitrary, and let $t \ge 4$ be such that $t-2$ is a prime power. Let $G$ be a graph with chromatic number $\chi(G) = t$. Then
\[
\gamma^+(S^{r}(G)) = \frac{(t-3)(t-2)^{r-1}}{(t-2)^r-1}.
\]
\end{theorem*}

Our proof is based on a double-counting argument that turns global positive codegree constraints into inequalities involving the neighborhood of a covered set. Then a weighted graph inequality yields the upper bound. The weighted graph inequality is sharp and might be of independent interest. The lower bound construction is the same as in \cite{da_co}, but we present it here nonetheless so that the article can be self-contained.

This result also shows the connection between daisy-free hypergraphs and matroid basis hypergraphs. When $t-2$ is a prime power, van der Pol, Walsh and Wigal \cite{Tmbh} proved that the construction of Ellis, Ivan, and Leader in \cite{da_co} is optimal within the class of matroid basis hypergraphs, with respect to Tur\'an density. Here we show that the same construction is optimal for positive codegree density among all $\mathcal D_{r,t}$-free hypergraphs, showing the connection between excluding a fixed matroid minor and daisy-free hypergraphs.

The article is organized as follows. In Section 2, we introduce the double-counting inequality and illustrate its use. In Section 3.1 we prove the upper bound for $\gamma^+(\mathcal D_{r,t})$ via a weighted graph inequality. In Section 3.2 we present the construction with a matching lower bound. Section 4 consists of some related open problems.

\section{Basic Tools and Examples}
In this section, $H$ shall be an $r$-uniform weighted hypergraph with $w_v>0$ denoting the weight of vertex $v$, and the total weight denoted $N$. For a vertex set $A$, define $w(A) = \prod_{v\in A}w_v$. For a covered set $F$ of size $r-1$, we define
\[
d(F)=\sum_{v\in N_H(F)}w_v
\]
where $N_H(F)$ is the neighborhood of $F$ in $H$. The corresponding positive codegree of $H$ shall be 
\[
\delta^+ = \min_{F \text{ covered}}d(F).
\]
We also define
\[
P = \sum_{e\in E(H)} w(e).
\]
For a covered set $F$ with $|F|=r-1$, fix $u\in F$ and consider the link graph of $F\setminus \{u\}$. Let $X=N_H(F)$. Then $u$ is adjacent to every vertex in $X$ in the link graph. Let $G_{F,u}$ denote the induced subgraph of the link graph on the vertex set $X$. We remark that if $H$ is $\mathcal D_{r,t}$-free, that is, every link graph is $K_t$-free, then every $G_{F,u}$ is $K_{t-1}$-free.
The $r-1$ graphs $G_{F,u}$ ($u\in F$) are all defined on the vertex set $X$. For $x,y\in X$, write $m(x,y)$ (or $m_F(x,y)$ if we want to specify $F$) for the number of graphs $G_{F,u}$ containing $xy$ as an edge. We set $m(x,x)=0$. We define the following link pairing sum
\[
L(F) = \sum_{x,y\in X}\frac{w_x w_y}{1+m(x,y)},
\]
which is one of the key elements in our proof. We emphasize that the sum is taken over ordered pairs $(x,y)$.
\begin{lemma}
\[
\sum_{|F|=r-1\text{ covered}}w(F)L(F)\leq NP.
\]
\end{lemma}
\begin{proof}
The proof is based on standard double counting. We know 
\[
\begin{aligned}
\sum_{F}w(F)L(F)&=\sum_F w(F)\sum_{x,y\in N_H(F)}\frac{w_x w_y}{1+m_F(x,y)}\\
&= \sum_{x}w_x\sum_{F\cup \{u\}\text{ is an edge},\,x\in N_H(F)}\frac{w(F\cup\{u\})}{1+m_F(u,x)}.
\end{aligned}
\]
We know $m_F(u,x) + 1 = \sum_{v\in F}\mathds{1}_{ux\in G_{F,v}}+1 = \sum_{v\in e(=F\cup \{u\})}\mathds 1 _{e\setminus \{v\}\cup \{x\}\in E(H)}$. Thus, 
\[\begin{aligned}
\sum_{F}w(F)L(F)& = \sum_{x}w_x\sum_{F\cup\{u\} = e\in E(H),\,x\in N_H(F)}\frac{w(e)}{\sum_{v\in e(=F\cup \{u\})}\mathds 1 _{e\setminus \{v\}\cup \{x\}\in E(H)}}\\&=\sum_{x}w_x\sum_{e\in E(H),\exists u\in e,\,x\in N_H(e\setminus \{u\})}w(e)\\
&\leq \sum_{x}w_x \sum_{e\in E(H)} w(e) = NP.
\end{aligned}\]
This proves the lemma.
\end{proof}
This lemma is the key point of our proof, and we shall investigate what we have lost in this inequality. The only inequality used in the proof is 
\[
\sum_{e\in E(H),\exists u\in e,\,x\in N_H(e\setminus \{u\})}w(e)\leq \sum_{e\in E(H)} w(e).
\]
So, if for every $e$ and $x\notin e$, there exists $u\in e$ such that $e\cup \{x\}\setminus \{u\}$ forms an edge, then we lose nothing in the inequality. Matroid basis hypergraphs have this property if the matroid is loopless.
\begin{corollary}\label{cor}
If 
\[
L(F)\geq \alpha d(F)^2
\]
for all covered $F$, where $\alpha>0$ is a constant independent of $F$, then 
\[
\delta^+(H)\leq \frac{N}{r\alpha}.
\]
\end{corollary}
\begin{proof}
We know that 
\[
\sum_{F}w(F)L(F)\leq NP
\]
and thus
\[
\sum_{F}\alpha w(F)d(F)^2\leq NP=N\sum_{F}\frac 1r w(F)d(F)
\]
where the equality follows by double counting. Since $d(F)\geq \delta^+$, we know
\[
\sum_{F}\alpha\delta^+ w(F)d(F)\leq N\sum_{F}\frac 1r w(F)d(F)
\]
and thus
\[
\delta^+\leq \frac N{r\alpha}.
\]
\end{proof}
This corollary allows us to focus on the local quantities $L(F)$, $d(F)$ to obtain a global bound on $\delta^+$. To see the power of these results, we give the following simple example. The result was already proved in \cite{pcdj} via a short structural argument, but our proof illustrates how the inequalities above can be used to study positive codegree density for graph suspensions.
\begin{theorem}
\[
\gamma^+(\mathcal D_{r,3})  = \frac1r.
\]
\end{theorem}
\begin{proof}
We know every link graph is triangle-free and thus the neighborhood of a vertex in the link graph must have no edges. Therefore, for a fixed $F$, each $m(x,y) = 0$ and thus $L(F) = \sum_{x,y\in N(F)} w_xw_y = d(F)^2$. So we conclude that $\gamma^+(\mathcal D_{r,3})\leq \frac 1r$. The complete balanced $r$-partite $r$-graph gives the matching lower bound.
\end{proof}

\section{Exact Value of
  \texorpdfstring{$\gamma^+(\mathcal D_{r,t})$}{gamma+(D(r,t))}}
\subsection{The Upper Bound}
The main theorem of this section is as follows.
\begin{theorem}\label{mth}
\[
\gamma^+(\mathcal D_{r,t})\le
 \frac{(t-3)(t-2)^{r-1}}{(t-2)^r-1}
\]
for all $t\ge 4$.
\end{theorem}
To prove this result, it suffices to prove the following result.
\begin{theorem}\label{auxthm}
Let $G_1,G_2,\dots,G_s$ be $s$ $K_{q+1}$-free graphs ($q\ge 2$) on the same weighted vertex set $V$ of total mass $W$. Use $m(x,y)$ to denote the number of graphs $G_i$ containing $xy$ as an edge. Then
\[
L_s = \sum_{u,v}\frac{w_uw_v}{1+m(u,v)}\ge b_s(q)W^2,
\]
where $b_s(q)=\frac{q-q^{-s}}{(q-1)(s+1)}$.
\end{theorem}
We will first use this theorem to get the general upper bound for positive codegree densities of daisies.
\begin{proof}[Proof of Theorem \ref{mth}]
For a fixed covered $(r-1)$-set $F$, we know that every $G_{F,u}$ is $K_{t-1}$-free, and thus (by taking $s = r-1$, $q  =t-2$):
\[
L(F) = L_{r-1}\geq b_{r-1}(t-2)d(F)^2.
\]
Thus, by Corollary \ref{cor}, 
\[
\gamma^+(\mathcal D_{r,t})\leq \frac 1 {rb_{r-1}(t-2)} = \frac{(t-3)(t-2)^{r-1}}{(t-2)^r-1}.
\]
\end{proof}
We will prove Theorem \ref{auxthm} by induction. We first cite a classical result regarding clique densities in graphs. The reader can refer to the classical result of Bollob\'as \cite{bcd} for a proof of the following theorem, or find the stronger result of Reiher in \cite{rcd}, which established the general clique density theorem. Since the result is standard, we omit the proof here.
\begin{theorem}\label{cdthm}
Let $G$ be a weighted graph with total weight $N$ and edge density $p = \frac 2 {N^2}\sum_{xy\in E(G)}w_x w_y$. Let 
\[
C_l = \sum_{H \text{ is an $l$-clique in } G}w(H)
\]
be the weighted number of $l$-cliques, where $w(H)$ is the product of the vertex weights in $H$. Then
\[
C_l\geq \big((l-1)p - (l-2)\big)\frac{N^l}{l^{l-1}}.
\]
\end{theorem}

The next key lemma will be used for the induction step. In particular, we shall replace the weights $w_u$ by $f_J(u)w_u$ to continue our induction.
\begin{lemma}\label{ind_l}
Let $G$ be a $K_{q+1} $-free weighted graph ($q\ge 2$) with total weight $1$. Let $p$ denote its edge density. For each $(q-1)$-clique $J$ of $G$, there is a function $f_J(u)\ge 0$ such that for $X_{uv}:=\sum_{J}f_J(u)f_J(v)$ the following inequalities hold:
\[
0\leq X_{uv}\le \mathds 1_{uv\notin E(G)},
\]
\[
T = \sum_{u,v\in V(G)}w_u w_v X_{uv}\geq (q-1)p-(q-2).
\]
If there are no cliques of size $q-1$, we set $X_{uv}=0$ for all $u,v$ and the above still holds.
\end{lemma}
\begin{proof}
If $G$ contains no $K_{q-1}$, then $(q-1)p-(q-2)\le 0$ and the results hold automatically.
Use $A_{u,J}$ for the indicator function that $J\subset N(u)$. Define
\[
D(u) = \sum_{J\text{ forms a $(q-1)$-clique}}w(J) A_{u,J},
\]
and
\[
f_J(u) = \frac{\sqrt{w(J)} A_{u,J}}{\sqrt{D(u)}},
\]
where 
\[
w(J) = \prod_{v\in J}w_v.
\]
We let $\frac 00 = 0$ in the following proof to simplify our explanation. So we let $f_J(u) = 0 $ when $D(u) = 0$. Since $G$ is $K_{q+1}$-free, for $uv\in E(G)$, the common neighbors of $u,v$ cannot contain $K_{q-1}$. Thus, $X_{uv} = 0$ for $uv\in E(G)$. We have
\[
\sum_{J}f_J(u)f_J(v)=\sum_Jw(J)A_{u,J}A_{v,J}\bigg/\bigg(\sqrt{D(u)D(v)}\bigg).
\]
We know 
\[
D(u)=\sum_Jw(J)A_{u,J}\ge \sum_Jw(J)A_{u,J}A_{v,J}
\]
and similarly
\[
D(v)\ge \sum_Jw(J)A_{u,J}A_{v,J}.
\]
Thus, 
\[
\bigg(\sqrt{D(u)D(v)}\bigg)\ge \sum_Jw(J)A_{u,J}A_{v,J}.
\]
So $X_{uv}\leq 1$ everywhere. Now the remaining task is to estimate $T$. We write 
$$Q(J)=\sum_{u}\frac{w_uA_{u,J}}{\sqrt{D(u)}} $$
and 
\[
T=\sum_{u,v,J}w_u w_v f_J(u)f_J(v)=\sum_{J}w(J)Q(J)^2.
\]
Now we define $\lambda(J) = \left( \prod_{j\in J}D(j)\right)^{\frac{1}{2(q-1)}}$. By Cauchy's inequality and the AM-GM inequality,
\[
\begin{aligned}
\sqrt{T\sum_{J}w(J)\lambda(J)^2}&\geq\sum_{J,u}w(J)w_u A_{u,J}\frac{\left( \prod_{j\in J}D(j)\right)^{\frac{1}{2(q-1)}}}{D(u)^{\frac12}}\\
&=\sum_{K\text{ a $q$-clique}}w(K)\prod_{k\in K}\left( D(k)\right)^{\frac{1}{2(q-1)}}\sum_{u\in K}\frac 1{D(u)^{\frac12+\frac1{2(q-1)}}}\\
&\ge \sum_{K\text{ is a $q$-clique}}qw(K).
\end{aligned}
\]
We also have
\[\begin{aligned}
\sum_{J}w(J)\lambda(J)^2 &= \sum_{J}\prod_{j\in J}w_jD(j)^{\frac 1{q-1}}.
\end{aligned}
\]
To estimate $\sum_{J}\prod_{j\in J}w_jD(j)^{\frac 1{q-1}}$, set $z_j = w_jD(j)^{\frac 1{q-1}}$. We use the following lemma to estimate the sum.
\begin{lemma}\label{wZi}(weighted Zykov's inequality; see \cite{lgmi})
Let $G$ be a $K_{q+1} $-free weighted graph ($q\ge 2$). Assume that $z_u\ge 0$ and $\sum_{u}z_u = Z$. Then
\[
\sum_{J \text{ is a $(q-1)$-clique}}\prod_{j\in J}z_j\leq q^{-(q-2)}Z^{q-1}.
\]
\end{lemma}
We note that the sum $\sum_{J}\prod_{j\in J}z_j$ is linear in $z_u,z_v$ if $uv\notin E$, so we can assume one of $z_u,z_v$ is zero when the sum attains its maximum. Repeating this operation, we can assume that $z$ is supported on a clique whose size is at most $ q$ since $G$ is $K_{q+1}$-free. Then
\[
\sum_{J}\prod_{j\in J}z_j\le \sum_{i=1}^{q}\prod_{j\neq i }{z_j}\le q^{-(q-2)}Z^{q-1}
\]
where the last inequality follows from Maclaurin's inequality.
Using Lemma \ref{wZi} and H\"older's inequality, we get
\[
\begin{aligned}
\sum_{J}w(J)\lambda(J)^2&\le q^{-(q-2)}\left(\sum_{u}w_uD(u)^{\frac1 {q-1}}\right)^{q-1}\\
&\le q^{-(q-2)}\left(\sum_{u}{w_u}\right)^{q-2}\left(\sum_{u}{w_uD(u)}\right)\\
&=q^{-(q-2)}\sum_{K\text{ is a $q$-clique}}qw(K).
\end{aligned}
\]
In the last equality we use standard double counting, and note that the total weight is $1$. Note that we have obtained 
\[
T\sum_{J}w(J)\lambda(J)^2\ge \left(\sum_{K\text{ is a $q$-clique}}qw(K)\right)^2.
\]
If $G$ has no $q$-clique, the inequality follows by Tur\'an's theorem. In fact, we can get $(q-1)p-(q-2)\le 0$. Thus, we finally obtain
\[
T\ge q^{q-1}\sum_{K\text{ is a $q$-clique}}w(K).
\]
By Theorem \ref{cdthm}, 
\[
\sum_{K\text{ is a $q$-clique}}w(K)\ge \left((q-1)p-(q-2)\right)q^{-(q-1)},
\]
and then
\[
T\ge (q-1)p-(q-2).
\]
\end{proof}
We now prove Theorem \ref{auxthm}.
\begin{proof}
By homogeneity, we assume that the total weight $W$ is $1$. We argue by induction on $s$. For $s = 1$, the result is just Tur\'an's theorem. Now assume that the result holds for $s-1$. For the induction step, we take away one graph $G_i$ and consider the sum
\[
l_i = \sum_{u,v}\frac{w_uw_v\mathds 1_{uv\notin E(G_i)}}{1+m(u,v)}.
\]
The $m(u,v)$ are determined by $G_j(j\neq i)$ when $uv\notin E(G_i)$. Applying Lemma \ref{ind_l} to $G_i$, we know
\[
X_{uv}\le \mathds 1_{uv\notin E(G_i)}
\]
and 
\[
T \ge (q-1)p-(q-2)
\]
where $p$ is the edge density of $G_i$. Then we apply the inductive hypothesis to
\[
\begin{aligned}
\sum_{u,v}\frac{w_uw_vX_{uv}}{1+m(u,v)} &= \sum_{J\text{ is a $(q-1)$-clique in $G_i$}}\sum_{u,v}\frac{w_uf_J(u)w_vf_J(v)}{1+m(u,v)}\\
&\ge\sum_{J}b_{s-1}(q)(\sum_{u}{w_uf_J(u)})^2\\
&=b_{s-1}(q)T.
\end{aligned}
\]
Here $T=\sum_{u,v}w_uw_vX_{uv}$. Then we can estimate
\[\begin{aligned}
l_i&\ge \sum_{u,v}\frac{w_uw_vX_{uv}}{1+m(u,v)}+\sum_{uv\notin E(G_i)}\frac{w_uw_v(1-X_{uv})}{s}\\
&\ge b_{s-1}(q)T+\frac{1-p-T}{s}\\
&=\frac{1-p}{s}+T\frac{1-q^{1-s}}{s(q-1)}\\
&\ge\frac{1-p}s + \frac{1-q^{1-s}}{s(q-1)}((q-1)p-(q-2))\\
&=(\frac{q-2}{q-1}-p)\frac{q^{1-s}}{s}+\frac 1{s(q-1)}.
\end{aligned}
\]
Since $p\le \frac{q-1}{q}$, we know
\[
l_i\ge \frac{1-q^{-s}}{s(q-1)}.
\]
We have
\[
\sum_{i=1}^{s}l_i = \sum_{u,v}w_uw_v\frac {s-m(u,v)}{1+m(u,v)} = (s+1)L_s-W^2.
\]
Thus
\[
L_s\ge \frac{1+s\frac{1-q^{-s}}{s(q-1)}}{s+1} = b_s(q)
\]
and then the theorem is proved.
\end{proof}
Now together with the argument at the beginning of this section, we get the general upper bound on the positive codegree density of daisies. We remark that the bound in Theorem \ref{auxthm} is tight for all $s,q$. In fact, let the vertex set be $V=[q]^s$ with total weight $1$ uniformly distributed on vertices. We know that each $v\in V$ is of the form $v = (v_i)_{1\le i\le s}$. Let $uv\in G_i$ iff $u_i\neq v_i$ for $1\le i\le s$. Then each $G_i$ is $q$-partite, thus $K_{q+1}$-free. And we can calculate $L_s = b_s(q)$. 

\subsection{Lower Bound Construction}
We directly use the construction in \cite{da_co}. For completeness, we briefly illustrate their construction based on finite fields and linear algebra. \par
Assume that $q=t-2$ is a prime power. Let $\mathbb F_q$ be the finite field of size $q$. We first construct a $(q^r-1)$-vertex hypergraph $H$ and then take its balanced blow-ups. Let $V(H)=\mathbb F_q^r\setminus \{0\}$. Let $v_1v_2\cdots v_r$ be in $E(H)$ iff the vectors $v_i(1\le i\le r)$ form a basis of $\mathbb F_q^r$. For a covered $(r-2)$-vertex set $v_1v_2\cdots v_{r-2}$, let $G$ be its link graph. Whether $u_1u_2$ is in $G$ only depends on their images in the quotient space $\mathbb F_q^r/\mathrm{span}(v_1,v_2,\cdots ,v_{r-2})$. So we can label the vertex set of $G$ by $\mathbb F_q^2$ with $uv\in E$ iff the images of $u,v$ in the quotient space are linearly independent. The vertices corresponding to the zero vector are isolated. We can partition $\mathbb F_q^2\setminus\{0\}$ into $q+1$ lines. Vertices in the same line are not adjacent, whereas vertices in different lines are adjacent. The number of vectors $v$ that extend the independent set $v_1,v_2,\cdots,v_{r-1}$ to a basis is $(q-1)q^{r-1}$, which is exactly the positive codegree of $H$. So for the balanced blow-ups, the normalized positive codegree is asymptotically $\frac{(t-3)(t-2)^{r-1}}{(t-2)^r-1}$.

We note that the above construction works for the suspensions of every graph $G$ satisfying $\chi(G) = t$, since every link graph of the above $H$ is $(t-1)$-partite. Consequently, we can now prove our main theorem, which we state again for the convenience of the reader.
\begin{theorem*}
Let $r \ge 3$ be arbitrary, and let $t \ge 4$ be such that $t-2$ is a prime power. Let $G$ be a graph with chromatic number $\chi(G) = t$. Then
\[
\gamma^+(S^{r}(G)) = \frac{(t-3)(t-2)^{r-1}}{(t-2)^r-1}.
\]
\end{theorem*}
\begin{proof}
The lower bound construction was given above. For $G=K_t$ the result follows immediately from the upper bound that we proved in Section $3.1$. It was proved in \cite{pco_b} that positive codegree density is invariant under blow-ups. Since $\chi(G) = t$, $S^r(G)$ is contained in the blow-up of $\mathcal D_{r,t}$. Thus, $\gamma^+(S^{r}(G))\le \gamma^+(S^{r}(K_t))$ and we get the general result.
\end{proof}

\section{Concluding Remarks and Open Questions}
We conclude the article with several related open problems. One natural question is to study the positive codegree density for all generalized daisies without restrictions on $t$.
\begin{question}
What is the positive codegree density of daisies $\mathcal D_{r,t}$ when $t\ge 8$ and $t-2$ is not a prime power? 
\end{question}
It would also be interesting to study the suspensions of hypergraphs, and to determine whether the methods treating suspensions of graphs can be extended to suspensions of hypergraphs.
\begin{question}
Can we determine $\gamma^+(S^r(G))$ when $G$ is a hypergraph?
\end{question}
Another long-standing open problem is to determine the exact Tur\'an density of daisies. 
\begin{question}
What is the Tur\'an density of $\mathcal D_r$?
\end{question}
In \cite{cj_dt}, Bollob\'as, Leader and Malvenuto conjectured that the Tur\'an density of $\mathcal D_3$ is $\frac12$. The lower bound is obtained by iterated blow-ups of the complement of the Fano plane. This conjecture remains open.

\section*{Acknowledgment}
The author is grateful to Sohail Farhangi for proofreading and giving many helpful suggestions.

\section*{AI disclosure}
The author worked with AI (specifically, OpenAI's model GPT-6.1-Sol) to obtain a proof regarding the case of ordinary daisies $\mathcal D_r$. Then the author worked independently to obtain a more general proof handling all generalized daisies $\mathcal D_{r,t}$. The proof presented here is the author's more general proof. The article is entirely written by the author, who takes full responsibility for all the content. AI was also used to correct typos and assist with English-language editing in the manuscript.

\bibliographystyle{abbrv}
\bibliography{ref}

\end{document}